\pdfoutput=1
\documentclass[11pt,a4paper,reqno]{amsart}

\usepackage[T1]{fontenc}
\usepackage[utf8]{inputenc}
\usepackage{lmodern}
\usepackage{amsmath,amssymb,amsthm,mathtools}
\usepackage{enumitem}
\usepackage[protrusion=true,expansion=false]{microtype}
\usepackage[a4paper,textwidth=150mm,textheight=235mm,
  centering,headheight=14pt,headsep=7mm,footskip=10mm]{geometry}
\usepackage[unicode,hidelinks]{hyperref}

\hypersetup{
  pdftitle={Plasticity of the closed unit ball of real c0},
  pdfauthor={Jiajun Tang},
  pdfsubject={Plastic metric spaces, nonexpansive bijections, and the closed unit ball of real c0}
}

\newcommand{\N}{\mathbb{N}}
\newcommand{\R}{\mathbb{R}}
\newcommand{\norm}[1]{\lVert#1\rVert}
\newcommand{\abs}[1]{\lvert#1\rvert}

\DeclareMathOperator{\supp}{supp}
\DeclareMathOperator{\sgn}{sgn}
\DeclareMathOperator{\id}{id}
\DeclareMathOperator{\med}{med}
\DeclareMathOperator{\ext}{ext}

\numberwithin{equation}{section}

\theoremstyle{plain}
\newtheorem{theorem}{Theorem}[section]
\newtheorem{lemma}[theorem]{Lemma}
\newtheorem{proposition}[theorem]{Proposition}
\newtheorem{corollary}[theorem]{Corollary}

\theoremstyle{definition}
\newtheorem{example}[theorem]{Example}
\newtheorem{remark}[theorem]{Remark}

\newlist{conditions}{enumerate}{1}
\setlist[conditions]{
  label=\textup{(\roman*)},ref=\roman*,
  leftmargin=2.8em,labelwidth=2.1em,labelsep=0.7em,
  align=left,topsep=5pt,itemsep=3pt,parsep=0pt,partopsep=0pt
}

\allowdisplaybreaks[1]
\title[Plasticity of the closed unit ball of real $c_0$]
  {Plasticity of the closed unit ball of real $c_0$}

\author{Jiajun Tang}
\address{Basic Course Department, Zhonghuan Information College Tianjin University of Technology, Tianjin 300380, China}
\email{tt20982025@163.com}
\thanks{The author was supported by the Scientific Research Seedling Project for Young Teachers of Zhonghuan Information College Tianjin University of Technology (2025KYLG03).}
\date{}


\begin{document}

\begin{abstract}
We prove that the closed unit ball $B$ of $c_0(\R)=C_0(\N,\R)$ is a plastic metric space. Every nonexpansive bijection $F\colon B\to B$ is the restriction of a surjective linear signed coordinate permutation of $c_0(\R)$, and in particular is an isometric homeomorphism. More generally, every nonexpansive surjection $T\colon B\to B$ sending each coordinate into the closed segment between zero and that coordinate, and injective on every finite coordinate cube, coincides with the identity map. Examples on $B$ show that neither surjectivity nor injectivity on finite cubes can be omitted.
\end{abstract}

\maketitle

\section{Introduction}

Let $(M,d)$ and $(N,\rho)$ be metric spaces. A mapping $F\colon M\to N$ is called \emph{nonexpansive} if $\rho(F(x),F(y))\le d(x,y)$ for all $x,y\in M$, \emph{noncontractive} if $\rho(F(x),F(y))\ge d(x,y)$ for all $x,y\in M$, and an \emph{isometry} if equality holds for every pair of points. Following Naimpally, Piotrowski, and Wingler~\cite{NPW06}, a metric space $(M,d)$ is called \emph{plastic} if every nonexpansive bijection $F\colon M\to M$ is an isometry. Passing to set-theoretic inverses shows that $(M,d)$ is plastic if and only if every noncontractive bijection $G\colon M\to M$ is an isometry.

In metric-space topology, plasticity goes back to the classical theorem of Freudenthal and Hurewicz~\cite{FH36}, which states that every nonexpansive surjection of a compact metric space onto itself is an isometry. Naimpally, Piotrowski, and Wingler~\cite{NPW06} extended plasticity to all totally bounded metric spaces. Without total boundedness, boundedness and completeness no longer exclude strict compression of distances: an infinite-dimensional Hilbert ellipsoid with suitably chosen semi-axes admits a nonexpansive bijection onto itself that fails to be an isometry~\cite[Example~2.7]{CKOW16}. Topologically, a continuous bijection of a non-compact complete metric space onto itself need not be open or closed, so its inverse may fail to be continuous.

Cascales, Kadets, Orihuela, and Wingler~\cite{CKOW16} posed the problem of whether the closed unit ball $B_X$ of every Banach space $X$ is plastic. By the Mazur--Ulam theorem~\cite{MU32} and its convex-body extension due to Mankiewicz~\cite{Man72}, every surjective isometry of $B_X$ onto itself extends to a surjective linear isometry on $X$. The unit-ball plasticity problem therefore asks when bijective $1$-Lipschitz self-maps of $B_X$ inherit linear isometric rigidity, in parallel with Tingley's problem for unit spheres~\cite{Tin87,KM12,THL13,Ban22} and the sphere-plasticity problem~\cite{LZ26,Zav26}. Positive results for $B_X$ have been obtained for strictly convex spaces~\cite{CKOW16}, for $\ell_1$~\cite{KZ16}, for $\ell_1$-sums of strictly convex spaces~\cite{KZ18}, for spaces whose unit sphere is covered by finite-dimensional polyhedral extreme subsets~\cite{AKZ19}, and for several related spaces~\cite{ZLH22}. In all of these settings, extreme points play a central role because $F^{-1}(\ext(B_X))\subseteq\ext(B_X)$ for every nonexpansive bijection $F\colon B_X\to B_X$~\cite[Theorem~2.3]{CKOW16}.

This mechanism vanishes when $\ext(B_X)=\varnothing$. Let $\N=\{0,1,2,\ldots\}$ carry the discrete topology, and let $\alpha\N=\N\cup\{\infty\}$ denote its Alexandroff one-point compactification. For the space $C(\alpha\N,\R)\cong c$ of convergent real sequences, the closed unit ball still has an abundance of extreme points, and Leo~\cite{Leo22} proved that $B_c$ is plastic; Fakhoury~\cite{Fak24} generalized this result to $C(K)$ when $K$ is a compact metrizable space with finitely many accumulation points. By contrast, the closed subspace
\[
  c_0 := c_0(\R) = C_0(\N,\R) \cong \{f\in C(\alpha\N,\R) : f(\infty)=0\}
\]
of real sequences tending to zero, equipped with the supremum norm
\[
  \norm{x} := \sup_{j\in\N}\abs{x_j},
\]
has closed unit ball $B:=\{x\in c_0:\norm{x}\le 1\}$ satisfying $\ext(B)=\varnothing$.

In~\cite[Section~3]{Leo22}, Leo established a signed coordinate normalization for nonexpansive bijections of $B$ onto itself and proved that $F\colon B\to B$ is an isometry provided that $F^{-1}$ is continuous, that is, provided that $F$ is a homeomorphism of $B$. Zavarzina~\cite{Zav26} recently proved that the unit sphere $S_{c_0}:=\{x\in c_0:\norm{x}=1\}$ is plastic, while noting that whether continuity of $F^{-1}$ can be omitted on $B$ remained open. Because a nonexpansive bijection $F\colon B\to B$ satisfies $F^{-1}(S_{c_0})\subseteq S_{c_0}$ but is not known beforehand to satisfy $F(S_{c_0})\subseteq S_{c_0}$, plasticity of the sphere does not transfer directly to the ball.

In the present paper we remove the continuity assumption on $F^{-1}$ and prove that every nonexpansive bijection of $B$ onto itself is the restriction of a signed coordinate permutation. Write $e_j$ for the $j$-th unit vector in $c_0$, $T_j(x)=(T(x))_j$ for the $j$-th coordinate of $T(x)$, and $\supp(x):=\{j\in\N:x_j\ne 0\}$ for the support of $x\in c_0$. For each finite subset $S\subset\N$, let
\[
  Q_S := \{x\in B : \supp(x)\subseteq S\},
\]
which is isometric to $([-1,1]^S,\norm{\cdot}_\infty)$ via extension by zero. We call a mapping $T\colon B\to B$ \emph{normalized} if
\begin{equation}\tag{N}\label{eq:segment}
  \min(0,x_j) \le T_j(x) \le \max(0,x_j)
  \qquad (x\in B,\ j\in\N).
\end{equation}
Condition~\eqref{eq:segment} implies $T(Q_S)\subseteq Q_S$ for every finite $S\subset\N$, yet allows a nonzero input coordinate to vanish in the image. In Section~2 we establish the unique signed coordinate normalization on $B$ together with a quantitative local openness property for continuous injections on finite cubes derived from Brouwer's fixed point theorem. In Section~3 we use two positive preimages and cubical openness to prove that every normalized nonexpansive surjection injective on each $Q_S$ preserves supports and therefore equals $\id_B$, deduce the plasticity of $B$ and its topological consequences, and give three examples testing the hypotheses.

\section{Coordinate normalization and finite cubes}

Throughout the paper, every nonzero sequence $x\in c_0$ attains its supremum norm at some index, since $\{j\in\N:\abs{x_j}\ge\norm{x}/2\}$ is finite and nonempty. We write $\sgn(t)=1$ for $t>0$, $\sgn(0)=0$, and $\sgn(t)=-1$ for $t<0$. For each $p\in B$, put
\[
  D(p) := \{z\in B : \norm{z-p}\le 1\}.
\]
For each $n\in\N$, the coordinate truncation $P_n\colon B\to Q_{\{0,\ldots,n-1\}}$ defined by $(P_n x)_j=x_j$ for $j<n$ and $(P_n x)_j=0$ for $j\ge n$ satisfies
\begin{equation}\label{eq:trunc-conv}
  \norm{x - P_n x} = \sup_{j\ge n}\abs{x_j} \longrightarrow 0
  \qquad (n\to\infty,\ x\in B).
\end{equation}
Thus $B_{00}:=\bigcup_{\text{finite }S\subset\N}Q_S$ is norm-dense in $B$.

We first record the properties of signed coordinate permutations and normalized mappings that will be used repeatedly.

\begin{lemma}\label{lem:basic}
The following statements hold:
\begin{conditions}
  \item For every bijection $\sigma\colon\N\to\N$ and every sequence $\alpha=(\alpha_i)_{i\in\N}\in\{-1,1\}^\N$, the formula
  \[
    (Ua)_{\sigma(i)} = \alpha_i a_i \qquad (a\in c_0,\ i\in\N)
  \]
  defines a surjective linear isometry $U\colon c_0\to c_0$ with inverse $(U^{-1}b)_i=\alpha_i b_{\sigma(i)}$, satisfying $U(B)=B$ and $U(Q_S)=Q_{\sigma(S)}$ for every finite $S\subset\N$.
  \item If $T\colon B\to B$ satisfies~\eqref{eq:segment}, then $x_j T_j(x)\ge 0$, $\abs{T_j(x)}\le\abs{x_j}$, $T(0)=0$, $\supp(T(x))\subseteq\supp(x)$, and $T(Q_S)\subseteq Q_S$ for every finite $S\subset\N$. Moreover, if $T_j(x)\ne 0$, then $\sgn(x_j)=\sgn(T_j(x))$, and $T_j(x)=s\in\{-1,1\}$ implies $x_j=s$.
  \item If $T\colon B\to B$ satisfies~\eqref{eq:segment} and $d=(d_j)_{j\in\N}\in\{-1,1\}^\N$, then the diagonal reflection $(D_d a)_j=d_j a_j$ is an isometric involution of $B$ preserving each $Q_S$, and $K=D_d\circ T\circ D_d$ also satisfies~\eqref{eq:segment}. Furthermore, $K$ is nonexpansive, surjective onto $B$, or injective on every finite cube $Q_S$ whenever $T$ has the corresponding property.
\end{conditions}
\end{lemma}

\begin{proof}
\textup{(i)} Given $a\in c_0$, define $y\in\R^\N$ by $y_k=\alpha_{\sigma^{-1}(k)}a_{\sigma^{-1}(k)}$. For every $\varepsilon>0$, the set $\{k\in\N:\abs{y_k}\ge\varepsilon\}=\sigma(\{i\in\N:\abs{a_i}\ge\varepsilon\})$ is finite because $\sigma$ is bijective and $a\in c_0$. Hence $y\in c_0$ and $\norm{y}=\sup_{i\in\N}\abs{\alpha_i a_i}=\norm{a}$. Linearity of $U\colon a\mapsto y$ is clear, and the same argument applied to $\sigma^{-1}$ shows that $(Vb)_i=\alpha_i b_{\sigma(i)}$ defines a linear isometry inverse to $U$, with $\supp(Ua)=\sigma(\supp(a))$.

\textup{(ii)} Fix $x\in B$ and $j\in\N$. Condition~\eqref{eq:segment} gives $T_j(x)=0$ when $x_j=0$, $0\le T_j(x)\le x_j$ when $x_j>0$, and $x_j\le T_j(x)\le 0$ when $x_j<0$. Thus $x_j T_j(x)\ge 0$ and $\abs{T_j(x)}\le\abs{x_j}$ in all cases, whence $\supp(T(x))\subseteq\supp(x)$ and $T(Q_S)\subseteq Q_S$. If $T_j(x)\ne 0$, then $x_j\ne 0$ and $x_j T_j(x)>0$, so $\sgn(x_j)=\sgn(T_j(x))$. If $T_j(x)=s\in\{-1,1\}$, then $1=\abs{T_j(x)}\le\abs{x_j}\le 1$ gives $x_j=s$.

\textup{(iii)} By \textup{(i)}, $D_d$ is a surjective isometry of $B$ satisfying $D_d^2=\id_B$ and $D_d(Q_S)=Q_S$, so $K=D_d\circ T\circ D_d$ inherits nonexpansiveness, surjectivity, and injectivity on each $Q_S$. For any $a\in B$ and $j\in\N$, condition~\eqref{eq:segment} at $D_d a$ gives $\min(0,d_j a_j)\le T_j(D_d a)\le\max(0,d_j a_j)$. Multiplying by $d_j\in\{-1,1\}$ maps the closed interval between $0$ and $d_j a_j$ onto the closed interval between $0$ and $a_j$, which is~\eqref{eq:segment} for $K$.
\end{proof}

Nonexpansive surjections fix the origin and carry coverings of $B$ by unit-radius balls $D(p)$ to coverings by the corresponding image balls; compare~\cite[Lemma~3.1]{Leo22}.

\begin{lemma}\label{lem:covers}
Let $F\colon B\to B$ be a nonexpansive surjection.
\begin{conditions}
  \item $F(0)=0$, and $B=\bigcup_{\lambda\in\Lambda}D(p_\lambda)$ implies $B=\bigcup_{\lambda\in\Lambda}D(F(p_\lambda))$ for every family $\{p_\lambda\}_{\lambda\in\Lambda}\subseteq B$.
  \item For nonzero $p,q\in B$, one has $B=D(p)\cup D(q)$ if and only if $\supp(p)=\supp(q)=\{i\}$ and $p_i q_i<0$ for some $i\in\N$.
  \item For $x\in B_{00}$ with $A=\supp(x)$ and $s_i=\sgn(x_i)$ ($i\in A$), one has
  \begin{equation}\label{eq:finite-cover}
    B = D(x) \cup \bigcup_{i\in A} D(-s_i e_i).
  \end{equation}
\end{conditions}
\end{lemma}

\begin{proof}
\textup{(i)} Given $y\in B$, choose $x\in B$ with $F(x)=y$; then $\norm{y-F(0)}\le\norm{x}\le 1$. If $F_k(0)\ne 0$ for some $k\in\N$, taking $y=-\sgn(F_k(0))e_k\in B$ would give
\[
  \norm{y - F(0)} \ge \abs{-\sgn(F_k(0)) - F_k(0)} = 1 + \abs{F_k(0)} > 1,
\]
a contradiction. Thus $F(0)=0$. Moreover, $F(D(p_\lambda))\subseteq D(F(p_\lambda))$ by nonexpansiveness, so surjectivity yields $B=F(B)\subseteq\bigcup_{\lambda\in\Lambda}D(F(p_\lambda))\subseteq B$.

\textup{(ii)} If $\supp(p)=\supp(q)=\{i\}$ and $p_i q_i<0$, then for each $z\in B$ we can choose $w\in\{p,q\}$ with $w_i z_i\ge 0$, which gives $\abs{z_i-w_i}\le 1$ and $\abs{z_j-w_j}=\abs{z_j}\le 1$ for $j\ne i$, so $z\in D(w)$. Conversely, suppose $B=D(p)\cup D(q)$ for nonzero $p,q\in B$. If $p_i\ne 0$ and $q_j\ne 0$ for distinct indices $i\ne j$, then $z=-\sgn(p_i)e_i-\sgn(q_j)e_j\in B$ satisfies $\norm{z-p}\ge 1+\abs{p_i}>1$ and $\norm{z-q}\ge 1+\abs{q_j}>1$, missing both balls. Since $\supp(p)$ and $\supp(q)$ are nonempty, both supports must equal the same singleton $\{i\}$. If $p_i q_i>0$, the vector $-\sgn(p_i)e_i\in B$ similarly misses both balls. Hence $p_i q_i<0$.

\textup{(iii)} Given $z\in B$, if $z_i x_i<0$ for some $i\in A$, then $\abs{z_i+s_i}=1-\abs{z_i}\le 1$ and $z\in D(-s_i e_i)$. Otherwise $z_i x_i\ge 0$ for all $i\in A$, so $\abs{z_j-x_j}\le 1$ for every $j\in\N$ and $z\in D(x)$.
\end{proof}

Once coordinate signs and zero coordinates are preserved, surjectivity forces each output coordinate to contract toward zero in modulus.

\begin{lemma}\label{lem:modulus-bound}
Let $H\colon B\to B$ be a nonexpansive surjection satisfying $x_j H_j(x)\ge 0$ for all $x\in B$, $j\in\N$, and $H_j(x)=0$ whenever $x_j=0$. Then $H$ satisfies~\eqref{eq:segment}.
\end{lemma}

\begin{proof}
Observe first that $H_j(u)\ne 0$ implies $u_j\ne 0$ and $u_j H_j(u)>0$, so $\sgn(u_j)=\sgn(H_j(u))$. Fix $x\in B$ and $i\in\N$. The case $x_i=0$ holds by hypothesis, so assume $x_i\ne 0$. Define $y^-,y^+\in B$ by $y_i^-=-1$, $y_i^+=1$, and for $j\ne i$,
\[
  y_j^- = y_j^+ =
  \begin{cases}
    (j+1)^{-1}, & x_j\ge 0,\\
    -(j+1)^{-1}, & x_j<0.
  \end{cases}
\]
By surjectivity, choose $u,v\in B$ with $H(u)=y^-$ and $H(v)=y^+$. For every $j\ne i$, the nonzero numbers $u_j$ and $v_j$ share the same strict sign, whence $\abs{u_j-v_j}<1$. Since
\[
  2 = \norm{y^- - y^+} = \norm{H(u) - H(v)} \le \norm{u - v} \le 2,
\]
we obtain $\abs{u_i-v_i}=2$. Because $u_i<0<v_i$, it follows that $u_i=-1$ and $v_i=1$.

Set $w=u$ if $x_i>0$, and $w=v$ if $x_i<0$. For $j\ne i$, either $x_j\ge 0<w_j\le 1$ or $-1\le w_j<0$ and $x_j<0$, so $\abs{x_j-w_j}\le 1$, whereas $w_i=H_i(w)=-\sgn(x_i)$. Consequently, $\norm{x-w}=1+\abs{x_i}$. Since $x_i H_i(x)\ge 0$, we also have $\abs{H_i(x)-H_i(w)}=1+\abs{H_i(x)}$, and nonexpansiveness yields
\[
  1 + \abs{H_i(x)} = \abs{H_i(x) - H_i(w)} \le \norm{x - w} = 1 + \abs{x_i}.
\]
Thus $\abs{H_i(x)}\le\abs{x_i}$, which together with $x_i H_i(x)\ge 0$ gives~\eqref{eq:segment}.
\end{proof}

Combining Lemmas~\ref{lem:basic}--\ref{lem:modulus-bound} yields the unique signed coordinate normalization of Leo~\cite[Section~3]{Leo22}.

\begin{proposition}\label{prop:normal}
For every nonexpansive bijection $F\colon B\to B$, there exist a unique bijection $\sigma\colon\N\to\N$, a unique sequence $\alpha=(\alpha_i)_{i\in\N}\in\{-1,1\}^\N$, and a normalized nonexpansive bijection $H\colon B\to B$ such that
\begin{equation}\label{eq:factor-F}
  F_{\sigma(i)}(x) = \alpha_i H_i(x)
  \qquad (x\in B,\ i\in\N).
\end{equation}
\end{proposition}

\begin{proof}
By Lemma~\ref{lem:covers}\textup{(i)} and injectivity of $F$, we have $F(0)=0$ and $F(p)\ne 0$ for $p\ne 0$. Fix $i\in\N$. Applying Lemma~\ref{lem:covers}\textup{(i)}--\textup{(ii)} first to the covering pair $\{e_i,-e_i\}$ and then to $\{te_i,-\sgn(t)e_i\}$ for $0<\abs{t}\le 1$, we obtain an index $\sigma(i)\in\N$, a sign $\alpha_i=\sgn(F_{\sigma(i)}(e_i))\in\{-1,1\}$, and a continuous function $g_i\colon[-1,1]\to[-1,1]$ such that
\[
  F(te_i) = g_i(t)e_{\sigma(i)} \quad (t\in[-1,1]),
  \qquad g_i(0)=0,
  \qquad \alpha_i t\,g_i(t)>0 \quad (t\ne 0).
\]
If $\sigma(i)=\sigma(j)=k$ for some $i\ne j$, then $g_i(\alpha_i)>0$ and $g_j(\alpha_j)>0$, so the intermediate value theorem gives $c>0$ and nonzero $t_1,t_2\in[-1,1]$ with $F(t_1 e_i)=ce_k=F(t_2 e_j)$, contradicting injectivity of $F$. Hence $\sigma\colon\N\to\N$ is injective.

Next, let $x\in B_{00}$ have finite support $A=\supp(x)$ and put $s_i=\sgn(x_i)$ for $i\in A$. By~\eqref{eq:finite-cover} and Lemma~\ref{lem:covers}\textup{(i)}, the sets $D(F(x))$ and $D(F(-s_i e_i))$ ($i\in A$) cover $B$. Define $y^{(0)}\in B$ by $y^{(0)}_{\sigma(i)}=\alpha_i s_i$ for $i\in A$ and $y^{(0)}_k=0$ for $k\notin\sigma(A)$. Since $F(-s_i e_i)=g_i(-s_i)e_{\sigma(i)}$ with $\sgn(g_i(-s_i))=-\alpha_i s_i\ne 0$, we have
\[
  \norm{y^{(0)} - F(-s_i e_i)} \ge \abs{\alpha_i s_i - g_i(-s_i)} = 1 + \abs{g_i(-s_i)} > 1
  \qquad (i\in A).
\]
Thus $y^{(0)}$ misses every $D(F(-s_i e_i))$ and must belong to $D(F(x))$, giving $\alpha_i s_i F_{\sigma(i)}(x)\ge 0$ for all $i\in A$. For any $k\notin\sigma(A)$, both $y^{(0)}+e_k$ and $y^{(0)}-e_k$ belong to $B$ and likewise miss every $D(F(-s_i e_i))$, so $\abs{\pm 1 - F_k(x)}\le 1$ and hence $F_k(x)=0$.

Now let $x\in B$ be arbitrary. By~\eqref{eq:trunc-conv} and continuity of $F$, $F(P_n x)\to F(x)$. If $x_i\ne 0$, then $(P_n x)_i=x_i$ for $n>i$, so passing to the limit gives $\alpha_i x_i F_{\sigma(i)}(x)\ge 0$. If $x_i=0$, injectivity of $\sigma$ gives $\sigma(i)\notin\sigma(\supp(P_n x))$ for all $n$, whence $F_{\sigma(i)}(x)=0$. Similarly, any $k\notin\sigma(\N)$ would satisfy $F_k(x)=0$ for all $x\in B$, contradicting surjectivity of $F$ at $e_k$; thus $\sigma$ is bijective. Taking $U$ from Lemma~\ref{lem:basic}\textup{(i)} and setting $H=U^{-1}\circ F\colon B\to B$, we see that $H$ is a nonexpansive bijection satisfying $x_i H_i(x)\ge 0$ and $H_i(x)=0$ when $x_i=0$, so $H$ satisfies~\eqref{eq:segment} by Lemma~\ref{lem:modulus-bound}.

Finally, if $(\widetilde{\sigma},\widetilde{\alpha},\widetilde{H})$ also satisfies~\eqref{eq:factor-F} with $\widetilde{H}$ normalized, then evaluating at $e_i$ and using $F(e_i)\ne 0$ and Lemma~\ref{lem:basic}\textup{(ii)} gives $F(e_i)=\widetilde{\alpha}_i\widetilde{H}_i(e_i)e_{\widetilde{\sigma}(i)}$ with $\widetilde{H}_i(e_i)>0$, which forces $\widetilde{\sigma}(i)=\sigma(i)$ and $\widetilde{\alpha}_i=\alpha_i$.
\end{proof}

We turn to continuous injections on a finite coordinate cube. For finite nonempty $S\subset\N$, write $Q=[-1,1]^S$ with norm $\norm{a}_\infty=\max_{j\in S}\abs{a_j}$, $Q^\circ=(-1,1)^S$, and $\partial Q=Q\setminus Q^\circ=\{b\in Q:\norm{b}_\infty=1\}$.

\begin{lemma}\label{lem:cubical-retraction}
Let $S$ be finite and nonempty, and put $Q=[-1,1]^S$. There exists no continuous mapping $R\colon Q\to\partial Q$ satisfying $R(b)=b$ for all $b\in\partial Q$.
\end{lemma}

\begin{proof}
Equip $\R^S$ with the Euclidean norm $\norm{a}_2=(\sum_{j\in S}a_j^2)^{1/2}$ and write $B_2=\{v\in\R^S:\norm{v}_2\le 1\}$. Define $\Phi\colon Q\to B_2$ and $\Psi\colon B_2\to Q$ by $\Phi(0)=\Psi(0)=0$ and
\[
  \Phi(a) = \frac{\norm{a}_\infty}{\norm{a}_2}\,a \quad (a\ne 0),
  \qquad
  \Psi(v) = \frac{\norm{v}_2}{\norm{v}_\infty}\,v \quad (v\ne 0).
\]
Then $\norm{\Phi(a)}_2=\norm{a}_\infty$ and $\norm{\Psi(v)}_\infty=\norm{v}_2$, so $\Phi$ and $\Psi$ are mutually inverse continuous mappings between $Q$ and $B_2$. If a continuous mapping $R\colon Q\to\partial Q$ satisfied $R(b)=b$ on $\partial Q$, then $\Phi\circ(-R)\circ\Psi\colon B_2\to B_2$ would have a fixed point $v_0\in B_2$ by Brouwer's fixed point theorem~\cite[Corollary~2.15]{Hat02}. Setting $a_0=\Psi(v_0)\in Q$ gives $a_0=-R(a_0)\in\partial Q$, whence $a_0=R(a_0)=-a_0$, contradicting $\norm{a_0}_\infty=1$.
\end{proof}

Lemma~\ref{lem:cubical-retraction} yields a quantitative ball-inclusion estimate and local openness for normalized continuous injections on $Q$.

\begin{proposition}\label{prop:cube}
Let $S$ be finite and nonempty, let $Q=[-1,1]^S$ carry $\norm{\cdot}_\infty$, and let $f\colon Q\to Q$ be a continuous injection satisfying $\min(0,a_j)\le f_j(a)\le\max(0,a_j)$ for all $a\in Q$ and $j\in S$. For $z\in Q^\circ$, put $\delta_z=\min_{b\in\partial Q}\norm{f(b)-f(z)}_\infty$. Then $\delta_z>0$ and
\begin{equation}\label{eq:ball-in-image}
  B_\infty(f(z),\delta_z) := \{y\in\R^S : \norm{y - f(z)}_\infty < \delta_z\} \subseteq f(Q^\circ).
\end{equation}
Moreover, for every neighborhood $V$ of $z$ in $Q$, the set $f(V)$ contains an open neighborhood of $f(z)$ in $\R^S$.
\end{proposition}

\begin{proof}
Compactness of $\partial Q\ne\varnothing$ and injectivity of $f$ ensure that $f(z)\notin f(\partial Q)$, so the minimum $\delta_z$ is attained and strictly positive. Note also that $f(0)=0$. Fix $y\in\R^S$ with $\norm{y-f(z)}_\infty<\delta_z$, and write $\rho(v)=v/\norm{v}_\infty$ for $v\in\R^S\setminus\{0\}$. On $\partial Q\times[0,1]$, consider the three continuous mappings
\begin{align*}
  P_1(b,s) &= (1-s)b + s f(b),\\
  P_2(b,s) &= f(b) - f(sz),\\
  P_3(b,s) &= f(b) - f(z) - s(y - f(z)).
\end{align*}
We verify that none of $P_1,P_2,P_3$ vanishes on $\partial Q\times[0,1]$. Given $b\in\partial Q$, choose $j\in S$ with $\abs{b_j}=1$. The interval condition gives $f_j(b)=\lambda b_j$ with $\lambda\in[0,1]$, so $\abs{(P_1(b,s))_j}=(1-s+s\lambda)\ge 1-s>0$ for $0\le s<1$, while $P_1(b,1)=f(b)\ne f(0)=0$ by injectivity. Next, $sz\in Q^\circ$ for all $s\in[0,1]$, so $sz\ne b$ and injectivity gives $P_2(b,s)\ne 0$. Finally, $\norm{P_3(b,s)}_\infty\ge\delta_z-s\norm{y-f(z)}_\infty>0$ for all $s\in[0,1]$. Since $P_1(b,1)=f(b)=P_2(b,0)$ and $P_2(b,1)=f(b)-f(z)=P_3(b,0)$, concatenating $\rho\circ P_1$, $\rho\circ P_2$, and $\rho\circ P_3$ on $[0,1]$, $[1,2]$, and $[2,3]$ defines a continuous homotopy $L\colon\partial Q\times[0,3]\to\partial Q$ with $L(b,0)=b$ and $L(b,3)=\rho(f(b)-y)$.

Suppose that $y\notin f(Q)$. Writing $r=\norm{a}_\infty$, define $R\colon Q\to\partial Q$ by
\[
  R(a) =
  \begin{cases}
    \rho(f(4a) - y), & 0\le r\le 1/4,\\
    L(a/r,\, 4(1-r)), & 1/4\le r\le 1.
  \end{cases}
\]
At $r=1/4$, we have $a/r=4a\in\partial Q$ and $4(1-r)=3$, so both formulas give $\rho(f(4a)-y)$. Thus $R\colon Q\to\partial Q$ is continuous and satisfies $R(b)=L(b,0)=b$ for all $b\in\partial Q$, contradicting Lemma~\ref{lem:cubical-retraction}. Hence $y\in f(Q)$, and $\norm{y-f(z)}_\infty<\delta_z$ excludes preimages in $\partial Q$, proving~\eqref{eq:ball-in-image}.

For the second assertion, choose a relatively open subset $W$ of $Q$ with $z\in W\subseteq V$. Compactness of $Q\setminus W$ and injectivity of $f$ imply that $f(Q\setminus W)$ is closed in $\R^S$ and misses $f(z)$. Removing $f(Q\setminus W)$ from $B_\infty(f(z),\delta_z)$ gives an open neighborhood of $f(z)$ in $\R^S$ contained in $f(W)\subseteq f(V)$.
\end{proof}

\section{Proof of the main results}

Proposition~\ref{prop:cube} allows us to rule out coordinate vanishing whenever two suitable strictly positive preimages are available.

\begin{lemma}\label{lem:two}
Let $T\colon B\to B$ be nonexpansive, satisfy~\eqref{eq:segment}, and be injective on $Q_S$ for every finite $S\subset\N$. Fix $i\in\N$, and suppose $u,v\in B$ have strictly positive coordinates and satisfy
\begin{equation}\label{eq:two-peaks}
  T_i(u) = 1,
  \qquad
  T_j(v) = 1 \quad \text{whenever } u_j = 1.
\end{equation}
If $x\in B$ has nonnegative coordinates and $x_i>0$, then $T_i(x)>0$.
\end{lemma}

\begin{proof}
Suppose toward a contradiction that $T_i(x)=0$, the only alternative permitted by~\eqref{eq:segment}. Since $x,u,v\in c_0$, we may choose a finite initial segment $S=\{0,1,\ldots,N\}$ containing $i$ such that $0\le x_j,u_j,v_j\le 1/2$ for all $j\in\N\setminus S$. Define $z\in Q_S$ by setting $z_j=0$ for $j\notin S$ and
\[
  m_j = \min\!\left\{x_j,\, u_j,\, v_j,\, \tfrac{1}{2}\right\},
  \qquad
  z_j = \frac{1 + m_j}{2}
  \qquad (j\in S).
\]
Then $1/2\le z_j\le 3/4$ on $S$, so $z\in Q_S^\circ$. If $j\in S$ and $x_j>0$, then $0<m_j\le 1/2$ and every $a\in\{x_j,u_j,v_j\}$ satisfies $m_j\le a\le 1$, which gives
\begin{equation}\label{eq:strict-half}
  \abs{z_j - a} \le \frac{1 - m_j}{2} < \frac{1}{2}.
\end{equation}
If $j\in S$ and $x_j=0$, then $z_j=1/2$ and $\abs{z_j-a}\le 1/2$ for each $a\in\{x_j,u_j,v_j\}$. Together with the tail bounds outside $S$, we obtain $\norm{z-x}\le 1/2$, $\norm{z-u}\le 1/2$, and $\norm{z-v}\le 1/2$. Comparing $T_i(z)\in[0,1]$ with $T_i(x)=0$ and $T_i(u)=1$ via nonexpansiveness forces $T_i(z)=1/2$.

Put $J=\{j\in S:\abs{z_j-u_j}=1/2\}$. For each $j\in J$, the strict inequality~\eqref{eq:strict-half} excludes $x_j>0$, so $x_j=0$ and $z_j=1/2$; positivity of $u_j\in(0,1]$ then forces $u_j=1$. In particular, $i\notin J$ because $x_i>0$. By~\eqref{eq:two-peaks}, $\norm{z-v}\le 1/2$, and~\eqref{eq:segment}, every $j\in J$ satisfies
\begin{equation}\label{eq:pin-J}
  \frac{1}{2} = T_j(v) - \frac{1}{2} \le T_j(z) \le z_j = \frac{1}{2},
\end{equation}
so $T_j(z)=1/2$ for all $j\in J$.

Identify $Q_S$ with $[-1,1]^S$ via zero-extension. By Lemma~\ref{lem:basic}\textup{(ii)} and hypothesis, $T|_{Q_S}\colon Q_S\to Q_S$ is a continuous injection satisfying the hypotheses of Proposition~\ref{prop:cube}. Since $S\setminus J$ is finite, the set
\[
  V = \left\{a\in Q_S : \abs{a_j - u_j} < \frac{1}{2}\ \text{for all } j\in S\setminus J\right\}
\]
is a relatively open neighborhood of $z$ in $Q_S$. By Proposition~\ref{prop:cube}, there exist $t>0$ and $w\in V$ such that $T(w)=T(z)-te_i$, where equality holds off $S$ because both images vanish outside $S$.

We now bound $\abs{u_j-w_j}$ for all $j\in\N$. If $j\in J$, then $i\notin J$ and~\eqref{eq:pin-J} give $T_j(w)=T_j(z)=1/2$, so Lemma~\ref{lem:basic}\textup{(ii)} yields $1/2\le w_j\le 1=u_j$ and hence $\abs{u_j-w_j}\le 1/2$. On $S\setminus J$, membership $w\in V$ gives $\abs{u_j-w_j}<1/2$. Outside $S$, we have $w_j=0$ and $0<u_j\le 1/2$, so $\abs{u_j-w_j}\le 1/2$ as well. Consequently, $\norm{u-w}\le 1/2$, whereas $T_i(u)=1$ and $T_i(w)=1/2-t$ yield
\[
  \frac{1}{2} + t = \abs{T_i(u) - T_i(w)} \le \norm{u - w} \le \frac{1}{2},
\]
contradicting $t>0$. Therefore $T_i(x)>0$.
\end{proof}

Once supports are preserved everywhere on $B$, backward iterates of any finitely supported point remain confined to a single finite cube and converge monotonically.

\begin{lemma}\label{lem:orbit-id}
Let $T\colon B\to B$ be a nonexpansive surjection satisfying~\eqref{eq:segment}. If $\supp(T(x))=\supp(x)$ for every $x\in B$, then $T=\id_B$.
\end{lemma}

\begin{proof}
Fix a finite subset $S\subset\N$ and $a\in Q_S$. The case $S=\varnothing$ is immediate since $T(0)=0$, so assume $S\ne\varnothing$. By surjectivity, set $b^{(0)}=a$ and choose $b^{(m+1)}\in B$ recursively with $T(b^{(m+1)})=b^{(m)}$ for each $m\in\N$. Support preservation keeps $\supp(b^{(m)})=\supp(a)\subseteq S$, so $b^{(m)}\in Q_S$ for all $m\in\N$.

For each fixed $j\in S$, write $c_m=b^{(m)}_j$. By~\eqref{eq:segment}, $c_m$ lies in the closed interval between $0$ and $c_{m+1}$. Hence $(\max(0,c_m))_{m\in\N}$ is nondecreasing in $[0,1]$ and $(\min(0,c_m))_{m\in\N}$ is nonincreasing in $[-1,0]$, so their sum $c_m$ converges in $[-1,1]$. Because $S$ is finite, $(b^{(m)})_{m\in\N}$ converges in norm to some $p\in Q_S$. Continuity of $T$ gives $T(p)=p$, and nonexpansiveness of $T^m$ yields
\[
  \norm{a - p} = \norm{T^m(b^{(m)}) - T^m(p)} \le \norm{b^{(m)} - p} \longrightarrow 0
  \qquad (m\to\infty).
\]
Thus $a=p$ and $T(a)=a$. Since $B_{00}$ is dense in $B$ by~\eqref{eq:trunc-conv}, continuity of $T$ implies $T=\id_B$.
\end{proof}

Lemmas~\ref{lem:two} and~\ref{lem:orbit-id} establish our rigidity theorem for normalized nonexpansive surjections.

\begin{theorem}\label{thm:rigid}
Let $T\colon B\to B$ be a nonexpansive surjection satisfying~\eqref{eq:segment}. If $T|_{Q_S}$ is injective for every finite subset $S\subset\N$, then $T=\id_B$.
\end{theorem}

\begin{proof}
Fix $x\in B$ and $i\in\N$ with $x_i\ne 0$. Define $d=(d_j)_{j\in\N}\in\{-1,1\}^\N$ by $d_j=-1$ when $x_j<0$ and $d_j=1$ when $x_j\ge 0$. By Lemma~\ref{lem:basic}\textup{(iii)}, $K=D_d\circ T\circ D_d\colon B\to B$ is a nonexpansive surjection satisfying~\eqref{eq:segment} and injective on each $Q_S$.

Define $q^{(1)}\in B$ by $q^{(1)}_i=1$ and $q^{(1)}_j=2^{-j-1}$ for $j\ne i$, and choose $u\in B$ with $K(u)=q^{(1)}$. Lemma~\ref{lem:basic}\textup{(ii)} gives $u_j\ge q^{(1)}_j>0$ for all $j\in\N$ and $u_i=1$. Since $u\in c_0$, the set $P=\{j\in\N:u_j=1\}$ is finite. Define $q^{(2)}\in B$ by setting $q^{(2)}_j=1$ for $j\in P$ and $q^{(2)}_j=2^{-j-1}$ for $j\notin P$, and choose $v\in B$ with $K(v)=q^{(2)}$. Then $v_j\ge q^{(2)}_j>0$ for all $j\in\N$, and $K_j(v)=1$ whenever $u_j=1$. Applying Lemma~\ref{lem:two} to $K$ at $D_d x$, whose coordinates are $\abs{x_j}$, yields
\[
  d_i T_i(x) = K_i(D_d x) > 0,
\]
so $T_i(x)\ne 0$. Together with Lemma~\ref{lem:basic}\textup{(ii)}, this gives $\supp(T(x))=\supp(x)$ for all $x\in B$, and Lemma~\ref{lem:orbit-id} concludes the proof.
\end{proof}

Taking the contrapositive of Theorem~\ref{thm:rigid} immediately locates a finite-support collision for any non-identity normalized surjection.

\begin{corollary}\label{cor:collision}
Let $T\colon B\to B$ be a nonexpansive surjection satisfying~\eqref{eq:segment}. If $T\ne\id_B$, then there exist a finite subset $S\subset\N$ and distinct points $x,y\in Q_S$ such that $T(x)=T(y)$.
\end{corollary}

\begin{proof}
If $T|_{Q_S}$ were injective for every finite $S\subset\N$, Theorem~\ref{thm:rigid} would force $T=\id_B$.
\end{proof}

Combining Proposition~\ref{prop:normal} with Theorem~\ref{thm:rigid} gives the main theorem on the closed unit ball of $c_0$.

\begin{theorem}\label{thm:main}
Let $F\colon B\to B$ be a nonexpansive bijection. There exist a unique bijection $\sigma\colon\N\to\N$ and a unique sequence $(\alpha_i)_{i\in\N}\in\{-1,1\}^\N$ such that
\begin{equation}\label{eq:main-formula}
  F_{\sigma(i)}(x) = \alpha_i x_i
  \qquad (x\in B,\ i\in\N).
\end{equation}
In particular, $F$ is the restriction to $B$ of the surjective linear isometry $U\colon c_0\to c_0$ given by $(Ux)_{\sigma(i)}=\alpha_i x_i$, and $\norm{F(x)-F(y)}=\norm{x-y}$ for all $x,y\in B$.
\end{theorem}

\begin{proof}
Let $(\sigma,\alpha)$ and $H=U^{-1}\circ F\colon B\to B$ be given by Proposition~\ref{prop:normal}. Since $H$ is a bijection of $B$ onto itself, it is surjective onto $B$ and injective on every $Q_S$. Theorem~\ref{thm:rigid} gives $H=\id_B$, so $F=U|_B$ satisfies~\eqref{eq:main-formula} and preserves distances by Lemma~\ref{lem:basic}\textup{(i)}. Uniqueness of $(\sigma,\alpha)$ follows from Proposition~\ref{prop:normal}.
\end{proof}

Theorem~\ref{thm:main} yields automatic homeomorphism for nonexpansive bijections, isometric rigidity for noncontractive bijections, and plasticity of closed balls in $C_0(\Gamma,\R)$ on countably infinite discrete spaces.

\begin{corollary}\label{cor:consequences}
Let $B$ be the closed unit ball of $c_0$.
\begin{conditions}
  \item Every nonexpansive bijection $F\colon B\to B$ is an isometric homeomorphism of $B$ onto itself satisfying $F(S_{c_0})=S_{c_0}$.
  \item Every noncontractive bijection $G\colon B\to B$ is the restriction of a surjective linear signed coordinate permutation of $c_0$, hence an isometry.
  \item For every countably infinite discrete space $\Gamma$ and every $R>0$, the closed ball $B_R(C_0(\Gamma,\R))=\{f\in C_0(\Gamma,\R):\sup_{\gamma\in\Gamma}\abs{f(\gamma)}\le R\}$ is plastic, and every nonexpansive bijection of $B_R(C_0(\Gamma,\R))$ onto itself is the restriction of a signed permutation of $\Gamma$.
\end{conditions}
\end{corollary}

\begin{proof}
\textup{(i)} By Theorem~\ref{thm:main}, $F=U|_B$ for a surjective linear isometry $U$ of $c_0$, so $F^{-1}=U^{-1}|_B$ is continuous and $F(S_{c_0})=S_{c_0}$.

\textup{(ii)} If $G\colon B\to B$ is a noncontractive bijection, its set-theoretic inverse $F=G^{-1}\colon B\to B$ is a nonexpansive bijection, so $F=U|_B$ by Theorem~\ref{thm:main} and $G=U^{-1}|_B$.

\textup{(iii)} Fix a bijection $\tau\colon\N\to\Gamma$ and define the linear isometry $\Psi\colon C_0(\Gamma,\R)\to c_0$ onto $c_0$ by $(\Psi f)_j=f(\tau(j))$. For any nonexpansive bijection $\widetilde{F}$ of $B_R(C_0(\Gamma,\R))$ onto itself, the map
\[
  F\colon B \longrightarrow B, \qquad F(x) = \frac{1}{R}\,\Psi\!\left(\widetilde{F}\!\left(R\,\Psi^{-1}(x)\right)\right),
\]
is a nonexpansive bijection. Applying Theorem~\ref{thm:main} to $F$ and using linearity of $\Psi$ and $U$ yields the representation on $B_R(C_0(\Gamma,\R))$.
\end{proof}

We now show that none of the main hypotheses in Theorem~\ref{thm:rigid} can be dropped. First, without surjectivity, even a globally injective normalized nonexpansive mapping on $B$ can erase unit coordinates on every signed coordinate face.

\begin{example}\label{ex:folding}
Define $P_+\colon[0,1]^2\to[0,1]^2$ by
\[
  P_+(p,q) =
  \begin{cases}
    \left(\tfrac{1}{3}p,\, \tfrac{1}{3}q\right), & p+q\le 1,\\[1ex]
    \left(\tfrac{1}{3}(1-q),\, \tfrac{1}{3}(p+2q-1)\right), & p+q\ge 1,
  \end{cases}
\]
and extend $P_+$ to $\widetilde{P}\colon[-1,1]^2\to[-1,1]^2$ via the $90^\circ$ rotation $\mathcal{R}(p,q)=(-q,p)$ by setting $\widetilde{P}(\mathcal{R}^r(p,q))=\mathcal{R}^r(P_+(p,q))$ for $(p,q)\in[0,1]^2$ and $r\in\{0,1,2,3\}$. Since both formulas for $P_+$ agree on $p+q=1$ and give $P_+(t,0)=(t/3,0)$ and $P_+(0,t)=(0,t/3)$ for $t\in[0,1]$, adjacent quadrants agree along all four half-axes. On $p+q\ge 1$ in $[0,1]^2$, we have $0\le\frac{1}{3}(1-q)\le\frac{1}{3}p\le p$ and $0\le\frac{1}{3}(p+2q-1)=\frac{1}{3}(2q-(1-p))\le\frac{2}{3}q\le q$, so $\widetilde{P}$ satisfies the coordinatewise segment condition and $\norm{\widetilde{P}(p,q)}_\infty\le 2/3$. On each of the eight closed triangles into which $p=0$, $q=0$, and $\abs{p}+\abs{q}=1$ divide $[-1,1]^2$, $\widetilde{P}$ is affine with Jacobian rows of $\ell_1$-norm at most $1$ (namely $1/3$ and $1/3+2/3=1$), so subdividing segments shows that $\widetilde{P}$ is $\norm{\cdot}_\infty$-nonexpansive.

On $[0,1]^2$, writing $(u,v)=P_+(p,q)$, we recover $(p,q)=(3u,3v)$ when $u+v\le 1/3$ and $(p,q)=(6u+3v-1,1-3u)$ when $u+v\ge 1/3$, so $P_+$ is injective. Moreover, $P_+([0,1]^2)$ meets the coordinate axes along $[0,1/3]\times\{0\}$, $\{0\}\times[0,1/3]$, and $\{0\}\times(1/3,2/3]$; under $\mathcal{R}^r$, adjacent quadrant images meet only on the four axis segments of length $1/3$ where their axis preimages coincide, so $\widetilde{P}$ is injective on $[-1,1]^2$. Defining $P\colon B\to B$ blockwise by $(P_{2m}(x),P_{2m+1}(x))=\widetilde{P}(x_{2m},x_{2m+1})$ ($m\in\N$) gives a globally injective nonexpansive mapping satisfying~\eqref{eq:segment} with $\sup_{x\in B}\norm{P(x)}=2/3$. Since $\widetilde{P}(1,1)=(0,2/3)$, $\widetilde{P}(-1,-1)=(0,-2/3)$, $\widetilde{P}(-1,1)=(-2/3,0)$, and $\widetilde{P}(1,-1)=(2/3,0)$, every signed face $\{x\in B:x_k=s\}$ ($k\in\N$, $s=\pm 1$) contains a point with $P_k(x)=0$.
\end{example}

Next, surjectivity onto $B$ in Theorem~\ref{thm:rigid} cannot be weakened to norm-dense range, even when support preservation and injectivity on every finite cube hold.

\begin{example}\label{ex:tail-median}
For $a\le c$ and $t\in\R$, write $\med(a,t,c)=\max(a,\min(t,c))$. Put $b_j(x)=\sup_{k>j}\abs{x_k}$ and $\varepsilon_j=2^{-j-1}$, and define $D\colon B\to B$ by
\[
  D_j(x) = \varepsilon_j x_j + (1 - \varepsilon_j)\med(-b_j(x),\, x_j,\, b_j(x))
  \qquad (x\in B,\ j\in\N).
\]
Because $\med(-b_j(x),x_j,b_j(x))$ lies between $0$ and $x_j$ and $\varepsilon_j>0$, the mapping $D$ satisfies~\eqref{eq:segment} and $\supp(D(x))=\supp(x)$ for all $x\in B$; since $\min$, $\max$, and $b_j$ are $1$-Lipschitz, $D$ is nonexpansive. If $x,y\in Q_S$ are distinct for a finite $S\subset\N$ and $j=\max\{k\in S:x_k\ne y_k\}$, then $b_j(x)=b_j(y)=:b$, and strict monotonicity of $t\mapsto\varepsilon_j t+(1-\varepsilon_j)\med(-b,t,b)$ gives $D_j(x)\ne D_j(y)$; hence $D|_{Q_S}$ is injective. For any $y\in Q_S$ and any $m\in\N$ greater than every element of $S$, the vector $x^{(m)}=y+\norm{y}e_m\in B$ satisfies $D(x^{(m)})=y+2^{-m-1}\norm{y}e_m\to y$ as $m\to\infty$, so $\overline{D(B)}=B$. Yet $e_0\notin D(B)$ because $D(x)=e_0$ would force $\supp(x)=\{0\}$ and $\norm{D(x)}\le\varepsilon_0=1/2$. Furthermore, setting $x_j=2^{-j-1}$, $d_0=1/32$, and $d_{j+1}=-\varepsilon_j(1-\varepsilon_j)^{-1}d_j$ ($j\in\N$), induction gives $\abs{d_j}\le\frac{1}{8}x_j$, so $x$ and $y=x+d$ are distinct strictly decreasing positive sequences in $B$. Thus $b_j(x)=x_{j+1}$ and $b_j(y)=y_{j+1}$, whence
\[
  D_j(y) - D_j(x) = \varepsilon_j d_j + (1 - \varepsilon_j)d_{j+1} = 0
  \qquad (j\in\N),
\]
showing that $D$ is not globally injective on $B$.
\end{example}

Conversely, injectivity on finite cubes cannot be omitted from Theorem~\ref{thm:rigid}: there exists a non-identity normalized nonexpansive surjection of $B$ onto itself, which necessarily exhibits the finite-cube collision of Corollary~\ref{cor:collision}.

\begin{example}\label{ex:record}
For $x\in B$ and $i\in\N$, put $L_i(x)=\max_{0\le k<i}\abs{x_k}$ (with $L_0(x)=0$), $R_i(x)=\sup_{k>i}\abs{x_k}$, and $t_+=\max(0,t)$, and define $C\colon B\to B$ by
\[
  C_i(x) = \sgn(x_i)\max\!\left\{\min(\abs{x_i},\, R_i(x)),\, \tfrac{1}{3}(\abs{x_i} - L_i(x))_+\right\}
  \qquad (x\in B,\ i\in\N).
\]
Clearly $C$ satisfies~\eqref{eq:segment}. For $x,y\in B$ and $i\in\N$, if $x_i y_i\le 0$, condition~\eqref{eq:segment} gives $\abs{C_i(x)-C_i(y)}\le\abs{x_i-y_i}\le\norm{x-y}$; if $x_i y_i>0$, oddness allows us to assume $x_i,y_i>0$, where $\min(x_i,R_i(x))$ is $1$-Lipschitz and $\frac{1}{3}(x_i-L_i(x))_+$ is $\frac{2}{3}$-Lipschitz, so $C$ is nonexpansive.

To see that $C(B)=B$, note that $C(0)=0$ and fix $y\in B\setminus\{0\}$. Inductively define right-record indices $r_1<r_2<\cdots$ and heights $b_1>b_2>\cdots>0$ by $b_1=\norm{y}$, $r_1=\max\{k\in\N:\abs{y_k}=b_1\}$, stopping at $M=m$ if $\sup_{k>r_m}\abs{y_k}=0$, and otherwise setting $b_{m+1}=\sup_{k>r_m}\abs{y_k}$ and $r_{m+1}=\max\{k>r_m:\abs{y_k}=b_{m+1}\}$. Define $x\in B$ to agree with $y$ except that $x_{r_{m+1}}=\sgn(y_{r_{m+1}})b_m$ for each valid $m\ge 1$, and $x_{r_M+1}=b_M$ when the sequence terminates at $r_M$. At non-record indices $i$ before the records end, we have $R_i(x)\ge\abs{x_i}$ and hence $C_i(x)=x_i=y_i$. At $i=r_1$, we have $R_{r_1}(x)=b_1=\abs{x_{r_1}}$ and $C_{r_1}(x)=y_{r_1}$. At $i=r_{m+1}$ ($m\ge 1$), we have $L_i(x)\ge\abs{x_{r_1}}=b_1\ge b_m=\abs{x_i}$, so $(\abs{x_i}-L_i(x))_+=0$, while $R_i(x)=b_{m+1}<b_m$ gives $C_i(x)=\sgn(y_i)b_{m+1}=y_i$. When records terminate at $r_M$, we similarly have $L_{r_M+1}(x)\ge b_1\ge b_M$ and $R_{r_M+1}(x)=0$, so $C_{r_M+1}(x)=0=y_{r_M+1}$, and $C_i(x)=0=y_i$ for $i>r_M+1$. Thus $C(x)=y$ and $C$ is surjective. Nevertheless, $C(te_0)=\frac{1}{3}te_0$ for $t\in[-1,1]$, and $C(t(e_0+e_1))=C(t(e_0+e_2))=te_0$ in $Q_{\{0,1,2\}}$ for every $t\in(0,1]$.
\end{example}

\begin{remark}\label{rem:open}
The proof of Theorem~\ref{thm:rigid} relies on $c_0$-decay twice: first to ensure that the peak set $P=\{j\in\N:u_j=1\}$ of the first positive preimage is finite, and second to confine the tail of $\{x,u,v\}$ within $[-1/2,1/2]$ outside a finite cube $Q_S$. Both properties fail in the non-separable space $\ell_\infty(\N)\cong C(\beta\N,\R)$ on the Stone--\v{C}ech compactification $\beta\N$, where plasticity of the closed unit ball remains open~\cite{Leo22}. More generally, in view of Fakhoury's theorem for $C(K)$~\cite{Fak24} and Corollary~\ref{cor:consequences}\textup{(iii)} for discrete spaces, it remains open to determine for which non-compact locally compact Hausdorff spaces $L$ the closed unit ball of $C_0(L,\R)$ is plastic without assuming continuity of the inverse.
\end{remark}

\end{document}